\documentclass[reqno,centertags, 12pt]{amsart}
\usepackage{amsmath,amsthm,amscd,amssymb}
\usepackage{latexsym,verbatim}

\newcommand{\bbR}{{\mathbb{R}}}

\newcommand{\bbP}{{\mathbb{P}}}
\newcommand{\bbE}{{\mathbb{E}}}

\newcommand{\calS}{{\mathcal S}}
\newcommand{\calM}{{\mathcal M}}

\newcommand{\lb}{\label}

\newcommand{\beq}{\begin{equation}}
\newcommand{\eeq}{\end{equation}}
\newcommand{\ba}{\begin{align}}
\newcommand{\ea}{\end{align}}
\newcommand{\eps}{\varepsilon}
\newcommand{\del}{\delta}
\newcommand{\tht}{\theta}
\newcommand{\ka}{\kappa}
\newcommand{\al}{\alpha}

\newcommand{\til}{\tilde}

\newcounter{smalllist}
\newenvironment{SL}{\begin{list}{{\rm\roman{smalllist})}}{%
\setlength{\topsep}{0mm}\setlength{\parsep}{0mm}\setlength{\itemsep}{0mm}%
\setlength{\labelwidth}{2em}\setlength{\leftmargin}{2em}\usecounter{smalllist}%
}}{\end{list}}

\DeclareMathOperator*{\dist}{dist}
\allowdisplaybreaks
\numberwithin{equation}{section}

\newtheorem{theorem}{Theorem}[section]

\newtheorem{proposition}[theorem]{Proposition}
\newtheorem{lemma}[theorem]{Lemma}
\newtheorem{corollary}[theorem]{Corollary}
\theoremstyle{definition}
\newtheorem{definition}[theorem]{Definition}

\theoremstyle{remark}

\begin{document}
\title{Quenching of combustion by shear flows}

\author{Alexander Kiselev }

\address{ Institute for Advanced Study \\ Princeton NJ 08540  and
Department of Mathematics \\ University of Wisconsin \\ Madison,
WI 53706, USA \\ email: kiselev@math.wisc.edu}

\author{Andrej Zlato\v s}

\address{Department of Mathematics \\ University of Wisconsin \\ Madison,
WI 53706, USA \\ \, email: zlatos@math.wisc.edu}


\maketitle

\begin{abstract}
We consider a model describing premixed combustion in the presence
of fluid flow: reaction diffusion equation with passive advection
and ignition type nonlinearity.
What kinds of velocity profiles are capable of quenching
(suppressing) any given flame, provided the velocity's amplitude
is adequately large? Even for shear flows, the solution turns out
to be surprisingly subtle.
In this paper, we provide a sharp characterization of quenching
for shear flows: the flow can quench any initial data if and only
if the velocity profile does not have an interval larger than a
certain critical size where it is identically constant.
The efficiency of quenching depends strongly on the geometry and
scaling of the flow. We discuss the cases of slowly and quickly
varying flows, proving rigorously scaling laws that have been
observed earlier in numerical experiments. The results require new
estimates on the behavior of the solutions to advection-enhanced
diffusion equation (also known as passive scalar in physical
literature), a classical model describing a wealth of phenomena in
nature. The technique involves probabilistic and PDE estimates, in
particular applications of Malliavin calculus and central limit
theorem for martingales.
\end{abstract}

\section{Introduction} \lb{S.1}

A mathematical model
that describes a chemical reaction in a fluid is a system of two
equations for concentration $n$ and temperature $T$ of the form
\begin{eqnarray}\label{gen-sys}
 && T_t+u\cdot\nabla T=\kappa\Delta T+Mg(T)n\\
&&  n_t+u\cdot\nabla n=\frac{\kappa}{\hbox{Le}}\Delta n-
Mg(T)n.\nonumber
\end{eqnarray}
Here $\kappa$ is the thermal diffusivity, Le the Lewis number (ratio
of thermal and material diffusivities), and $M$ the reaction
strength. The equations (\ref{gen-sys}) are coupled to the reactive
Navier-Stokes equations for the advection velocity $u(x,y,t)$.  Two
assumptions are usually made to simplify the problem: the first is a
constant density approximation \cite{CW} that allows to decouple the
Navier-Stokes equations from the system (\ref{gen-sys}) and to
consider $u(x,y,t)$ as a prescribed quantity that does not depend on
$T$ and $n$.  The second assumption  is that $\hbox{Le}=1$  (equal
thermal and material diffusivities). These two assumptions reduce
the above system to a single scalar equation for the temperature
$T$. We assume in addition that the advecting flow is
unidirectional. Then the system (\ref{gen-sys}) becomes
\begin{eqnarray}
  \label{eq:1.1}
  &&T_t+Au(y)T_x=\kappa\Delta T+Mf(T)\\
&& T(0,x,y)=T_0(x,y)\nonumber
\end{eqnarray}
with $f(T)=g(T)(1-T)$. We are interested in strong advection, and
accordingly have written the velocity as a product of the
amplitude $A$ and the profile $u(y)$. In this paper we consider a
nonlinearity $f\not\equiv 0$ of the ignition type
\begin{eqnarray}
\nonumber
&&\hbox{(i) }~~\hbox{$f(0)=f(1)=0$ and $f(T)$ is Lipschitz continuous on $[0,1]$},\\
&&\hbox{(ii) }~\exists\theta_0>0 \hbox{ such that $f(T)=0$ for
$T\in[0,\theta_0]$, $f(T)\ge 0$ for
$T \in (\theta_0,1),$ }\label{eq:2.1.2}\\
&&\hbox{(iii) }~~f(T)\le T.\nonumber
\end{eqnarray}
The last condition in (\ref{eq:2.1.2}) is just a normalization. We
consider the reaction-diffusion equation (\ref{eq:1.1}) in the
strip $D=\left\{x\in {\mathbb R} ,~ y\in [0,h]\right\}$. We take
$u(y)$ to be periodic with period $h$ and with mean equal to zero:
\begin{equation}
  \label{eq:3.1.1}
  \int_{0}^hu(y)dy=0.
\end{equation}
A constant non-zero mean can be easily taken into account by
translation.
For the temperature, we impose periodic boundary conditions
\begin{equation}
  \label{eq:2.1.1}
  T(t,x,y)=T(t,x,y+h)
\end{equation}
in $y$ and decay in $x.$ We will always assume that initial data
$T_0(x,y)$ is such that $0\le T_0(x,y)\le 1$. Then we have $0\le
T\le 1$ for all $t>0$ and $(x,y)\in D$. For simplicity, we will
usually assume that the initial data coincide with characteristic
function of some set. More generally, we may assume that for some
$L$ and $\eta>0$ we have
\begin{eqnarray}\label{eq:2.1.4}
&& T_0(x,y)>\theta_0+\eta~~\hbox{for $|x|\le L/2$},\\
&& T_0(x,y)=0~~\hbox{for $|x|\ge L$}. \nonumber
\end{eqnarray}

Equation (\ref{eq:1.1}) may be considered as a simple model of
flame propagation in a fluid \cite{Ber-Lar-Lions}, advected by a
shear (unidirectional) flow. The physical literature on the
subject is vast, and we refer to the recent review \cite{jxin-3}
for an extensive bibliography. The main physical effect of
advection on front-like solutions is the speed-up of the flame
propagation due to the large scale distortion of the front. The
role of the advection term in (\ref{eq:1.1}) for the front-like
initial data was also a subject of intensive mathematical scrutiny
recently, see \cite{Berrev, jxin-3} for the references.

The present paper considers advection effects for a different
physically interesting situation, where initial data are compactly
supported. In this case, two generic scenarios are possible. If
the support of the initial data is large enough, then two fronts
form and propagate in opposite directions. Fluid advection speeds
up the propagation, accelerating the burning. However, if the
support of the initial data is small, then the advection  exposes
the initial hot region to diffusion which cools it below the
ignition temperature $\theta_0$, ultimately extinguishing the
flame. The main purpose of this paper is to study the possibility
of quenching of flames by strong fluid advection in the model
(\ref{eq:1.1}). The phenomena associated with flame quenching are
of great interest for physical, astrophysical and engineering
applications. For example, modelling of quenching and propagation
of reaction fronts in fluid flow are relevant to studies of
internal combustion engines, nuclear burning in stars and forest
fires. Mathematically, the problem reduces to studying the
advection-enhanced dissipation rate for the passive scalar. The
passive scalar equation is one of the most studied PDE models, and
has been the subject of extensive research by both physicists and
mathematicians. However, the question that we address here -
controlling the rate of decay of the $L^\infty$ norm in terms of
the amplitude and geometric properties of the flow - while
extremely natural, remained largely open until some very recent
work.

The problem of extinction and flame propagation in the
mathematical model (\ref{eq:1.1})
was first studied by Kanel
\cite{Kanel'} in one dimension and with no advection. He showed that,
in the absence of fluid motion, there exist two length scales
$L_0<L_1$ such that the flame becomes extinct for $L<L_0$, and
propagates for $L>L_1$.  More precisely, he has shown that there exist
$L_0$ and $L_1$ such that
\begin{eqnarray}
  \label{eq:2.2}
&&  T(t,x,y)\to 0~\hbox{as $t\to\infty$ uniformly in $D$ if $L<L_0$}\\
&&T(t,x,y)\to 1~\hbox{as $t\to\infty$ for all $(x,y)\in D$ if $L>L_1$}.
\nonumber
\end{eqnarray}
We note in passing that it has been only very recently established
by one of the authors that $L_0=L_1$ in this situation
\cite{Zlatos2}. In the absence of advection, the flame extinction
is achieved by diffusion alone, given that the support of initial
data is small compared to the scale of the laminar front width $l
= \sqrt{\kappa/M}.$ However, in many applications quenching is the
result of a strong wind, intense fluid motion, and operates on
larger scales.
Kanel's result was extended to non-zero advection by shear flows
by Roquejoffre \cite{Roq-2} who has shown that (\ref{eq:2.2})
holds also for $u\ne 0$ with $L_0$ and $L_1$ depending, in
particular, on $A$ and $u(y)$. However the interesting question
about more explicit quantitative dependence of $L_0,$ $L_1$ on $A$
and $u(y)$ remained completely open until recent work \cite{CKR}.
The following definition was given in \cite{CKR}.

\begin{definition} \lb{D.1.1}
  We say that the profile $u(y)$ is {\it quenching} if for any $L$ and any
  initial data $T_0(x,y)$ supported inside the interval
  $[-L,L]\times[0,h],$ there exists $A_0=A_0(M,\kappa,f,u,L)$ such that the solution of
  (\ref{eq:1.1}) becomes extinct:
\begin{equation} \lb{1.2}
 T(t,x,y)\to 0~\hbox{as $t\to\infty$ uniformly in $D$}
\end{equation}
for all $|A|\ge A_0$. We call the profile $u(y)$ {\it strongly
quenching} if the critical amplitude of advection $A_0$ satisfies
$A_0 \leq C L$ for some constant $C=C(M,\kappa,f,u)$ (which has the
dimension of inverse time).
\end{definition}

The quenching property has been linked in \cite{CKR} to
hypoellipticity of a certain degenerate diffusion equation. In
particular, one of the main results showed that $u(y)$ is strongly
quenching if there is no point $y$ where all derivatives of $u$
vanish. On the other hand, if $u(y)$ has a plateau larger than a
certain critical size, then $u$ is not quenching. However
hypoellipticity does not provide a precise solution of the problem
at hand: a shear flow $u(y)$ with a small plateau leads to an
auxiliary equation which is not hypoelliptic, yet it is quenching.
The first main result of this paper, Theorem \ref{T.3.1}, provides
a  sharp characterization of quenching shear flows. It states that
a shear flow is quenching if and only if it has a plateau
exceeding certain critical size. Nearby plateaux of smaller size
will not lead to the same effect. This critical scale can be
described in terms of existence of solutions to a nonlinear
Dirichlet problem. The main new technical ingredient involves
estimates on certain stochastic integrals, in particular
application of Malliavin calculus to derive absolute continuity of
the relevant random variables.

The second goal is to study the dependence of quenching on scaling
of the flow. Numerical experiments \cite{VCKRR} suggest that there
is a certain scale of the flow for which quenching is most
efficient. Namely, if $u(y) = \sin \alpha y,$ then the size $L_A$
of initial data that can be quenched by flow $Au(y)$ satisfies
$L_A \sim C_\alpha A$ with $C_\alpha$ achieving maximum for some
$\alpha_0.$ Moreover, the constant $C_\alpha$ satisfies $C_\alpha
\sim \alpha^{-1}$ for large $\alpha$ and $C_\alpha \sim \alpha^2$
for small $\alpha.$ Our Theorems \ref{T.4.1} and \ref{T.4.2},
which apply to general shear flows, prove that in the small and
large $\alpha$ asymptotic regimes one indeed has quenching for the
initial data satisfying the above scaling. Central limit-type
theorem for martingales is instrumental in obtaining the large
$\alpha$ result.

We mention that in a separate work \cite{Zla} one of us
investigates the phenomenon of quenching in the presence of
combustion-type reaction functions that do not have an ignition
cutoff $\tht_0$, but are allowed to be positive for all $T\in
(0,1)$. An important example of such function is the Arrhenius
type reaction $f(T)= e^{-A/T}(1-T),$ which is used in modeling of
many chemical reactions. One of the main results of that paper,
Theorem 1.3, is related to our Theorem \ref{T.3.1}. It states that
the quenching property of shear flows is in this case linked not
only to the size of their plateaux but also to the decay rate of
$f$ at $T=0$. Namely, if (in our setting) $f(T)\ge cT^p$ with
$p<3$, then no flow is quenching, whereas if $f(T)\le cT^p$ with
$p>3$, then flows with small enough plateaux are quenching and
those with large plateaux are not. This result, however, does not
provide a sharp characterization of quenching flows, and should be
understood as an extension to non-ignition reactions of the
results of \cite{CKR} rather than of our Theorem \ref{T.3.1}.

As a final remark we note that proving results for the system
\eqref{gen-sys} is typically much harder than for a single
equation \eqref{eq:1.1}, due to the lack of appropriate comparison
principles. This is not so in our case. While in the paper we
discuss quenching for a single equation \eqref{eq:1.1}, all our
quenching results (including Theorems \ref{T.3.1}(i), \ref{T.4.1},
and \ref{T.4.2}) extend immediately to the case of the system
\eqref{gen-sys}. This is a consequence of a remnant of the maximum
principle: the concentration $n(t,x)$ remains bounded above by one
for all times. Then $T_t+u\cdot \nabla T \leq \kappa \Delta T +M
g(T)$ in \eqref{gen-sys}, and all the bounds we prove for
quenching in the single equation model apply.

The paper is organized as follows. In Section 2 we establish some
auxiliary technical estimates on stochastic integrals. In Section
3 we prove results on quenching by shear flows and provide a
characterization of the critical plateau size in terms of the
corresponding Dirichlet problem. In Section 4 we deal with the
scaling question.

\section{Stochastic Integrals} \lb{S.2}

Results from this section will be used to obtain upper bounds on
the solutions of \eqref{eq:1.1} without the non-linear term, which
can be expressed in terms of the Brownian motion. See the
beginning of Section \ref{S.3} for details and how this translates
into estimates on the temperature $T$.

We call a {\it plateau} of a function $u\in C(\bbR)$ any maximal
(w.r.t.~inclusion) interval on which $u$ is constant. We start by
proving

\begin{lemma} \lb{L.2.1}
Let $u\in C^1(\bbR)$ be bounded along with its first derivative
and let $W_s^y$ denote the normalized one-dimensional Brownian
motion starting at $y$. Then for any $a\in\bbR$ we have
\begin{equation} \lb{2.1}
\bbP \bigg( \int_0^t u(W_s^y)ds = a \bigg) = \bbP \bigg( u(W_s^y)
= \frac at \text{ for } s\in[0,t] \bigg).
\end{equation}
\end{lemma}

\smallskip
{\it Remarks.} 1. In other words, the first probability is zero
unless $y$ is an interior point of a plateau of $u$ with
$u(y)=\tfrac at$, in which case it equals the probability of
$W_s^y$ staying inside this plateau for all $s\in[0,t]$.

\smallskip
2. This lemma for $u\in C^\infty$ and $y$ not in a plateau of $u$
follows from a probabilistic version of H\" ormander's theorem
(see, e.g., \cite[Theorem 2.3.2]{Nu}). Here we extend it to all
$u\in C^1$ and all $y$.

\smallskip
3. We believe that the same result holds for $u\in C(\bbR)$ but we
were unable to locate an appropriate reference in the literature.

\begin{proof}
By Theorem~2.1.3 in \cite{Nu} with $F(W^y)\equiv\int_0^t
u(W^y_s)ds$, the law of the random variable $F$ is absolutely
continuous with respect to the Lebesgue measure on $\bbR$ whenever
\begin{equation} \lb{2.1a}
\big( \|DF\|_2^2 = \big) \int_0^t \bigg( \int_s^t u'(W^y_r)dr
\bigg)^2 ds
>0
\end{equation}
almost surely. We note that with the notation of
\cite[p.24-26]{Nu}, if $u\in C^1$, then $F\in\mathbb{D}^{1,1}$ is
the limit of $F_n(W^y)\equiv \tfrac 1n \sum_{k=1}^n
u(W^y_{tk/n})$, and $DF(s)=\int_s^t u'(W^y_r)dr$ is the limit of
\[
DF_n(s)=\frac 1n \sum_{k=1}^n u'(W^y_{tk/n})
\chi_{[0,\frac{tk}n)}(s) = \sum_{k=1}^n \Big[ \frac 1n\sum_{j=k}^n
u'(W^y_{tj/n}) \Big] \chi_{[\frac{t(k-1)}n,\frac{tk}n)}(s).
\]

Eq.~\eqref{2.1a} is obviously true if $u'$ is not identically zero
on an interval around $y$, that is, when $y$ is not inside a
plateau. In particular, for such $y$ and all $a$,
\begin{equation} \lb{2.2}
\bbP \bigg( \int_0^t u(W_s^y)ds = a \bigg) = 0.
\end{equation}

Now assume $y$ to be inside a plateau $I$. For any open interval
$J$ with rational end points not intersecting any plateau of $u$,
and any rational $\tau\in(0,t)$, let $B_{J,\tau}$ be the set of
Brownian paths $W^y$ such that $W_\tau^y\in J$. Notice that every
$W^y$ that exits $I$ before time $t$, belongs to some such
$B_{J,\tau}$.

We have for any $a$
\begin{equation}\label{2.2a}
\bbP \bigg( \int_0^t u(W_s^y)ds = a \,\bigg | \, W^y\in B_{J,\tau}
\bigg) = 0.
\end{equation}
This follows from \eqref{2.2} applied to the $\int_\tau^t$ portion
of the integral. Indeed --- since $U_s\equiv W_{s+\tau}^y$ (for
$s\ge 0$) is just Brownian motion starting at $W_\tau^y$, given
any history $\{W_s^y\}_{s\le\tau}$, the probability of
$\int_\tau^t u(W_s^y)ds \, (=\int_0^{t-\tau} u(U_s)ds)$ being
$a-\int_0^\tau u(W_s^y)ds$ is zero because $W_\tau^y$ is not in a
plateau of $u$ if $W^y\in B_{J,\tau}$. By Fubini's theorem,
\eqref{2.2a} holds. Since there are only countably many sets
$B_{J,\tau}$, the result follows.
\end{proof}

The main result of this section is

\begin{lemma} \lb{L.2.2}
Let $u\in C^1(\bbR)$ be periodic. Then for any compact interval
$S\subset (0,\infty)$ we have
\begin{equation} \lb{2.3}
\bbP \bigg( \int_0^t u(W_s^y)ds \in [a,a+\eps] \smallsetminus
\{tu(y)\} \bigg) \to 0
\end{equation}
as $\eps\to 0$, uniformly in $(t,y,a)\in S\times\bbR\times\bbR$.
\end{lemma}

\smallskip
{\it Remarks.} 1. Note that non-uniform convergence is an obvious
consequence of Lemma \ref{L.2.1}.

\smallskip
2. The importance of this lemma lies in the fact that for large
$A$ it gives us a uniform (in $(t,y,x)\in S\times\bbR\times\bbR$)
estimate on the solution of \eqref{3.3},\eqref{3.4} below, using
\eqref{3.6aa}.
Through \eqref{3.5},\eqref{3.6} this translates into an upper
bound on the temperature $T$.
\smallskip

To prove the lemma, consider the function
\[
p(t,y,a,\eps)\equiv \bbP \bigg( \int_0^t u(W_s^y)ds \in [a,a+\eps]
\smallsetminus \{ tu(y) \} \bigg),
\]
that is, the probability of $\int_0^t u(W_s^y)ds \in [a,a+\eps]$
and $\{u(W_s^y)\}_{s\le t}$ not constant.

\begin{lemma} \lb{L.2.3}
Under the conditions of Lemma \ref{L.2.2}, $p$ is jointly
continuous in $\bbR^+\times\bbR^2\times\bbR_0^+$.
\end{lemma}

\begin{proof}
For $\delta_1\in\bbR_0^+$ and $\delta_2,\delta_3,\delta_4\in\bbR$
let
\[
\delta\equiv \|u\|_\infty |\delta_1| + t\|u'\|_\infty |\delta_2| +
|\delta_3| + |\delta_4|.
\]
Then
\[
\bigg| \int_0^{t+\delta_1} u(W_s^y+\delta_2)ds - \int_0^t
u(W_s^y)ds \bigg| \le \|u\|_\infty |\delta_1| + t\|u'\|_\infty
|\delta_2|
\]
and we have
\begin{align*}
|p(t+ & \delta_1,
y+\delta_2,a+\delta_3,\eps+\delta_4)-p(t,y,a,\eps)|
\\ \le & \bbP \bigg ( \int_0^t u(W_s^y)ds\in
[a-\delta,a+\delta] \cup [a+\eps-\delta,a+\eps+\delta]
\smallsetminus \{ tu(y)\} \bigg )
\\ + & \bbP \big ( \text{exactly one of $\{u(W_s^y)\}_{s\le t}$
and $\{u(W_s^y+\delta_2)\}_{s\le t+\delta_1}$ is constant} \big).
\end{align*}
As $\delta\to 0$, the first probability goes to zero because by
Lemma \ref{L.2.1},
\[
\bbP \bigg ( \int_0^t u(W_s^y)ds\in \{a,a+\eps\} \smallsetminus \{
tu(y)\} \bigg ) = 0.
\]
The second probability goes to zero because
\[
\bbP \big ( \text{$\{u(W_s^y)\}_{s\le t}$ is constant} \big )
\]
is continuous in $(t,y)$.
\end{proof}

\begin{proof}[Proof of Lemma \ref{L.2.2}]
By Lemma \ref{L.2.1}, $p(t,y,a,0)=0$. Hence by Lemma \ref{L.2.3},
$p(t,y,a,\eps)\downarrow 0$ as $\eps\to 0$, for any $(t,y,a)$. By
joint continuity of $p$ we then have $p(t,y,a,\eps)\downarrow 0$
as $\eps\to 0$, uniformly in $(t,y,a)\in K$, for any compact
$K\subset\bbR^+\times\bbR^2$. But $p$ is periodic in $y$ and
$p(t,y,a,\eps)=0$ for $|a|>t\|u\|_\infty+\eps$. Thus
$p(t,y,a,\eps)\downarrow 0$ as $\eps\to 0$, uniformly in
$(t,y,a)\in S\times\bbR\times\bbR$, for any compact
$S\subset\bbR^+$.
\end{proof}

\section{The Quenching Flows} \lb{S.3}

Let $u(y)\in C^1(\bbR)$ be a periodic function and let $f(T)$ be
an ignition-type non-linearity satisfying (i)-(iii) of
\eqref{eq:2.1.2}. Let $T(t,x,y)$, $\Phi(t,x,y)$, and $\Psi(t,x,y)$
be the solutions of
\begin{align}
T_t & = \ka \triangle T - Au(y)T_x + M f(T) \label{3.1}
\\ \Phi_t & = \ka \triangle \Phi - Au(y) \Phi_x   \label{3.2}
\\ \Psi_t & = \ka \Psi_{yy} - Au(y)\Psi_x \label{3.3}
\end{align}
with $(t,x,y)\in \bbR^+_0\times\bbR^2$ and initial
conditions
\begin{equation}\label{3.4}
T(0,x,y) = \Phi(0,x,y) = \Psi(0,x,y) = \chi_{[-L,L]}(x).
\end{equation}

Notice that to prove quenching, one only needs to show
\begin{equation} \lb{3.7}
\|T(\tau,\cdot,\cdot)\|_\infty\le\theta_0
\end{equation}
for some $\tau>0$. Indeed, the maximum principle then implies
$T(t,x,y)\le\theta_0$ for all $t\ge\tau$. Hence we have
\begin{equation}\label{psT}
T_t = \ka \triangle T - Au(y)T_x
\end{equation}
for $t\ge\tau$. To show \eqref{1.2} we first notice that by
integrating \eqref{3.1} in $(x,y)\in\bbR\times[0,h]$ (where $h$ is
the period of $u$) we have for $\|\cdot\|_p\equiv
\|\cdot\|_{L^p(\bbR\times[0,h])}$
\[
\frac \partial{\partial t} \|T(t,\cdot,\cdot)\|_1 = M\int f(T) dxdy
\le M \|T(t,\cdot,\cdot)\|_1
\]
by \eqref{eq:2.1.2}, and so $\|T(\tau,\cdot,\cdot)\|_1<\infty$.
One can then, for instance, use the estimates on the parabolic
kernel of the operator $\triangle - u\cdot \nabla$ for periodic
divergence-free flow $u$ on $\bbR^n$ from \cite{No} to show that
the kernel for $\triangle - u\cdot \nabla$ on $\bbR\times[0,h]$
with periodic boundary conditions is bounded above by $Ct^{-1/2}$
for some $C$ and all $t>1$ (see \cite{Zla}). Therefore
\begin{equation}\label{dbT}
\|T(\tau+t,\cdot,\cdot)\|_\infty \le Ct^{-1/2}
\|T(\tau,\cdot,\cdot)\|_1
\end{equation}
for $t>1$ and \eqref{1.2} follows. Alternatively, there is a more
elementary proof of \eqref{dbT} based on proving a Nash-type
inequality for the evolution of \eqref{psT}, namely that
\[ \|T(\tau+t,\cdot,\cdot)\|_2 \le \tilde{C}t^{-1/2}
\|T(\tau,\cdot,\cdot)\|_1 \] with $\tilde{C}$ independent of the
flow. Such an estimate also leads to \eqref{dbT} by a duality
argument. See \cite{FKR,Nash} for more details.

The functions $\Phi$, $\Psi$ can be used to estimate the
non-linear evolution:
\begin{align}
T(t,x,y) & \le \Phi (t,x,y)e^{M t}    \label{3.5}
\\ \sup_x\Phi(t,x,y) & \le \sup_x\Psi(t,x,y).   \label{3.6}
\end{align}
The first bound is achieved by replacing $f(T)$ with $T$ in
\eqref{3.1}, while the second bound follows from the equality
\[
\Phi(t,x,y) = \int_{-\infty}^\infty G(t,x-x')\Psi(t,x',y)\,dx'
\]
where
\[
G(x,t)=\frac 1{\sqrt {4\pi\ka t}}\, e^{-x^2/4\ka t}
\]
is the fundamental solution of the one-dimensional heat equation
(using that $\|G(t,\cdot)\|_{L^1_x}=1$). The equality is verified by
plugging it into \eqref{3.2}.

Since $\Phi$ and $\Psi$ satisfy the above linear equations, we can
apply the results from the previous section to obtain the
following estimates. Let $(W^x,W^y)$ be the normalized
2-dimensional Brownian motion starting at $(x,y)$ and let
$(X^x_t,Y^y_t)$ be the random process starting at $(x,y)$ and
given by
\begin{align*}
dX^x_t & = \sqrt{2\ka} \, dW^x_t - Au(Y^y_t)dt,
\\ dY^y_t & = \sqrt{2\ka} \, dW^y_t.
\end{align*}
Thus, $Y^y_t=y+\sqrt{2\kappa} (W^y_t-y) = W^y_{2\kappa t}$ and
\[
X^x_t = x+\sqrt{2\ka} (W^x_t-x) - \int_0^t Au(Y^y_s)ds = W^x_{2\ka
t} - \frac{A}{2\ka}\int_0^{2\ka t} u(W^y_s)ds.
\]
Then we have by \eqref{3.2}, \eqref{3.4}, and Lemma 7.8 in
\cite{Oks},
\begin{equation*}
\Phi(t,x,y)  = \bbE \big( \Phi(0,X^x_t,Y^y_t) \big)
 = \bbP \bigg(
W_{2\ka t}^x-\frac{A}{2\ka}\int_0^{2\ka t} u(W_s^y)ds \in[-L, L]
\bigg).
\end{equation*}
Similarly,
\begin{equation} \label{3.6aa}
\Psi(t,x,y)  =  \bbP \bigg( x-\frac{A}{2\ka}\int_0^{2\ka t}
u(W_s^y)ds \in[-L, L] \bigg).
\end{equation}

The following result provides a sharp characterization of the
quenching flows (see Definition~\ref{D.1.1}).

\begin{theorem} \lb{T.3.1}
With the above notation, there exists $0<\ell< \infty$, depending
only on $M$, $\ka$, and $f$, such that the following hold.
\begin{SL}
\item[{\rm{(i)}}] If the longest plateau of $u$ is shorter than
$\ell$, then $u$ is quenching.
\item[{\rm{(ii)}}]
If the longest plateau of $u$ is longer than $\ell$, then $u$ is not
quenching.
\end{SL}
Moreover, this $\ell$ is the infimum of all $l$ such that the
equation
\begin{equation}\label{3.8}
\phi_t = \ka \triangle \phi + M f(\phi)
\end{equation}
on $(x,y)\in\bbR\times [0,l]$ with Dirichlet boundary conditions
at $y=0,l$, has a solution $\phi$ with $\phi(0,\cdot,\cdot)$
compactly supported (and taking values in $[0,1]$) such that
$\phi$ does not go uniformly to zero as $t\to\infty$.
\end{theorem}

The key step in the proof is the following proposition.

\begin{proposition} \lb{T.3.1a}
For any $l,L\ge 0$ let $\tau(l,L)$ be the minimal time such that any
solution $\phi$ of \eqref{3.8} on $(x,y)\in\bbR\times [0,l]$ with
Dirichlet boundary conditions at $y=0,l$ and $\phi(0,\cdot,\cdot)$
supported in $[-L,L]\times[0,l]$ (and taking values in $[0,1]$),
satisfies $\phi(t,x,y)\le \tht_0/2$ for $t\ge\tau(l,L)$. If such a
time does not exist, we set $\tau(l,L)=\infty$. Then with the above
notation and $l$ the length of the longest plateau of $u$ we have
the following.
\begin{SL}
\item[{\rm{(i)}}] If $\tau(l,L)<\infty$ for every $L<\infty$, then
$u$ is quenching.
\item[{\rm{(ii)}}] If $\tau(l,L_0)=\infty$ for some $L_0<\infty$, then
$u$ is not quenching (and for all $L\ge L_0$ and any $A$ the
temperature $T(t,x,y)$ does not go uniformly to zero as
$t\to\infty$).
\end{SL}
\end{proposition}

{\it Remark.} Note that this result also applies in the case
$l=\ell$. Whether quenching happens in this case depends not only on
whether solutions of \eqref{3.8} with Dirichlet boundary conditions
at $y=0,\ell$, initially compactly supported, go uniformly to zero,
but on this decay being uniform in all $\phi(0,\cdot,\cdot)$
supported in $[-L,L]\times[0,\ell]$ (for each $L$).

\begin{proof}
(ii) Without loss of generality we can assume that the longest
plateau of $u$ is $I=[0,l]$. Also without loss of generality, let
$u(0)=0$. Indeed --- if $u(0)\neq 0$ and $\til T$ is the solution of
\eqref{3.1} with $u(y)$ replaced by $\til u(y)=u(y)-u(0)$, then
$T(t,x,y)=\til T(t,x-Au(0)t,y)$ and the result for $\til T$
translates directly to $T$.

Assume that for $L=L_0$ the temperature $T$ (with initial condition
\eqref{3.4} and some $A$) goes uniformly to zero and let
$\tau_0<\infty$ be such that $\|T(\tau_0,\cdot,\cdot)\|_\infty
\le\tht_0/2$. Since by comparison theorems any $\phi$ from the
statement of the proposition must satisfy $\phi(x,y,t) \le T(x,y,t)$
for all $t$, we obtain $\tau(l,L_0)\le\tau_0$, a contradiction.
Therefore $T$ does not go uniformly to zero for $L\ge L_0$ and any
$A$, and so $u$ is not quenching.

(i) Fix $L<\infty$, choose any $0<\del<\min\{L,1\}$, and let
$\tau_2\equiv 1+\tau(l,L+\del)$. We will show that for large
enough $|A|$ and all $x,y$ one has
\begin{equation}\label{3.122}
T(\tau_2,x,y) \le\tht_0.
\end{equation}
This is \eqref{3.7} and so \eqref{1.2} will follow.

Our strategy to show \eqref{3.122} will be to estimate $T$ by $\Phi$
via Lemma~\ref{L.2.2} for $y$ outside of the plateaux of $u$ (with
large enough $|A|$), and by a suitable $\phi$ from the statement of
the proposition for $y$ inside a plateau.

Let $d$ be the Lipschitz constant for $f$, define $c\equiv M\max\{
d,1\}$, and pick $\tau_1\in(0,1)$ such that
\begin{equation}\label{3.11}
\bbP \bigg( |W_{2\ka\tau_1}^0|\ge \frac \del 2 \bigg) \le
\frac{\tht_0}{4e^{c\tau_2}}.
\end{equation}
This choice will become clear later.

Let $A_0$ be large enough so that for all $|A|\ge A_0$,
$t\in[\tau_1,\tau_2]$, $x\in\bbR$, and $y$ not in the interior of a
plateau of $u$ we have
\begin{equation} \lb{3.133}
\Psi(t,x,y)=\bbP \bigg( x-\frac{A}{2\ka}\int_0^{2\ka t} u(W_s^y)ds
\in[-L, L] \bigg) \le \frac{\tht_0}{2e^{2c\tau_2}}.
\end{equation}
This is possible by Lemma \ref{L.2.2} and the first remark after
Lemma \ref{L.2.1}. Then by \eqref{3.5} and \eqref{3.6}, for
$t\in[\tau_1,\tau_2]$ and $y$ not in the interior of a plateau,
\begin{equation}\label{3.15}
\sup_x T(t,x,y) \le e^{M t} \sup_x \Psi(t,x,y) \le e^{M t}
\frac{\tht_0}{2e^{2c\tau_2}} \le \frac{\tht_0}{2e^{c\tau_2}}.
\end{equation}
In particular, \eqref{3.122} holds when $y$ is outside of all
plateaux of $u$.

We are left with the case of $y$ inside a plateau. Hence consider
a plateau $I$ of maximal length $l$. All the following arguments
will also apply to any other plateau of length $\til l\le l$
because $\tau(\til l,L)\le \tau(l,L)$ by comparison theorems.
Therefore proving \eqref{3.122} for $y\in I$ will yield the same
statement for all other plateaux of $u$ and the proof will be
finished.

Again assume without loss of generality that $I=[0,l]$ and $u(0)=0$.
Increase $A_0$ (if necessary) so that for any $|A|\ge A_0$,
\begin{equation} \lb{3.144}
\bbP \bigg( x-\frac{A}{2\ka}\int_0^{2\ka\tau_1} u(W_s^y)ds
\in\bigg[-L-\frac \del 2, L+\frac \del 2\bigg] \bigg) \le
\frac{\tht_0}{4e^{c\tau_2}}
\end{equation}
whenever $y\in I$ and $|x|\ge L+\del$. Such $A_0$ exists because by
Lemma \ref{L.2.2},
\[
\bbP \bigg( \int_0^{2\ka\tau_1} u(W_s^y)ds \in \frac {2\ka}{A}
\bigg[-L-\frac \del 2+x, L+\frac \del 2+x\bigg] \bigg) \to 0
\]
as $|A|\to\infty$, uniformly in $y\in I$ and $x\notin [-L-\frac \del
2, L+\frac \del 2]$ (because $u(0)=0$). Using \eqref{3.5},
\eqref{3.11}, and \eqref{3.144}, it follows that for $y\in I$ and
$|x|\ge L+\del$,
\begin{align}
T(\tau_1,x,y) & \le e^{M\tau_1} \Phi(\tau_1,x,y)   \notag
\\ & = e^{M\tau_1} \bbP \bigg(
W_{2\ka\tau_1}^x-\frac{A}{2\ka}\int_0^{2\ka\tau_1} u(W_s^y)ds
\in[-L, L] \bigg) \notag
\\ & \le e^{M\tau_1} \bbP \bigg( x-\frac{A}{2\ka}\int_0^{2\ka\tau_1} u(W_s^y)ds \in\bigg[-L-\frac \del
2, L+\frac \del 2\bigg] \bigg) + e^{M\tau_1}
\frac{\tht_0}{4e^{c\tau_2}}  \notag
\\ & \le \frac{\tht_0}{2e^{c(\tau_2-\tau_1)}}    \label{3.14}
\end{align}

Next consider a function $\phi (t,x,y)$ defined on
$[\tau_1,\infty)\times\bbR\times I,$ taking values in $[0,1]$, and
satisfying \eqref{3.8} with Dirichlet boundary conditions at $y=0,l$
and initial data
\begin{equation}\label{3.13}
T(\tau_1,x,y) -  \frac{\tht_0}{2e^{c(\tau_2-\tau_1)}} \le \phi
(\tau_1,x,y) \le \chi_{[-L-\delta,L+\delta]}(x)\chi_I(y).
\end{equation}
Such a $\phi$ exists because of \eqref{3.15} and \eqref{3.14}, and
by the definition of $\tau_2$ we have
\begin{equation}\label{3.12}
\phi(\tau_2,x,y)\le \frac{\tht_0}2.
\end{equation}

Now let $\omega\equiv T-\phi$ for $(t,x,y)\in
[\tau_1,\infty)\times\bbR\times I$. Then
\begin{equation}\label{3.16}
\omega_t =\triangle\omega + M[f(T)-f(\phi)] \le \triangle\omega +
c|\omega|
\end{equation}
By \eqref{3.13},
\begin{equation}\label{3.17}
\omega(\tau_1,x,y)\le \frac{\tht_0}{2e^{c(\tau_2-\tau_1)}}
\end{equation}
and by \eqref{3.15},
\begin{equation}\label{3.18}
\sup_{x}\{\omega(t,x,0),\omega(t,x,l)\}\le
\frac{\tht_0}{2e^{c(\tau_2-\tau_1)}}
\end{equation}
for $t\in[\tau_1,\tau_2]$. Now for $\til\omega\equiv e^{-ct}\omega$
we have
\[
\til\omega_t \le \triangle\til\omega + c(|\til\omega|-\til\omega).
\]
Thus by \eqref{3.17}, \eqref{3.18}, and the maximum principle,
\[
\til\omega(t,x,y)\le \frac{\tht_0}{2e^{c\tau_2}},
\]
for $t\in[\tau_1,\tau_2]$. Thus $\omega(\tau_2,x,y)\le \tht_0/2$
whenever $y\in I$. So by \eqref{3.12}, $T(\tau_2,x,y)\le \tht_0$ for
$y\in I$. As mentioned before, this also holds for any other plateau
of $u$. Together with \eqref{3.15} this gives \eqref{3.122}, and the
result follows.
\end{proof}

\begin{proof}[Proof of Theorem~\ref{T.3.1}]
Let $\ell$ be defined as in the statement of Theorem~\ref{T.3.1}.
The fact that $\ell<\infty$ is proved in \cite{CKR} by
constructing a subsolution of \eqref{3.8} on $\bbR\times [0,l]$
for large enough $l$. Proposition \ref{P.3.3} below shows that
$\ell>0$. Notice that by comparison theorems (see e.g.
\cite[Chapter 10]{Sm}), a solution $\phi$ described in the
statement of Theorem \ref{T.3.1} exists when $l>\ell$ and does not
exist when $l<\ell$.

Let $l$ be the length of the longest plateau of $u$ and let
$\tau(l,L)$ be as in Proposition \ref{T.3.1a}.

(ii) If $l>\ell$, then there exists a solution $\phi$ described
above. This means that $\|\phi(t,\cdot,\cdot)\|_\infty > \tht_0$
for all $t$. Since $\phi(0,x,y)$ is supported in $[-L_0,L_0]\times
[0,l]$ for some $L_0<\infty$, we obtain $\tau(l,L_0)=\infty$.
Proposition \ref{T.3.1a}(ii) then gives (ii).

(i) If $l<\ell$, then take $\del\equiv(\ell-l)/3$. For any
$L<\infty$ let $\phi$ be a solution of \eqref{3.8} on $\bbR\times
[-\del,l+\del]$ with Dirichlet boundary conditions at
$y=-\del,l+\del$ and
\[
\chi_{[-L,L]}(x)\chi_{[0,l]}(y) \le \phi(0,x,y)\le
\chi_{[-L-1,L+1]}(x)\chi_{[-\del,l+\del]}(y).
\]
Since $l+2\del<\ell$, there is $\tau_L<\infty$ such that
$\|\phi(\tau_L,\cdot,\cdot)\|_\infty\le\tht_0/2$. By comparison
theorems this implies $\tau(l,L)\le\tau_L<\infty$. The result
follows from Proposition \ref{T.3.1a}(i).
\end{proof}

In \cite{CKR} an upper bound on $\ell$ was provided by
constructing a non-zero compactly supported $\phi(x,y)$ such that
\[
\ka \triangle \phi + M f(\phi) \ge 0
\]
in the sense of distributions. By comparison theorems, $\ell$ is
at most the diameter of the support of $\phi$. Here we give a
lower bound on $\ell$, in terms of the existence of a stationary
1D solution of \eqref{3.8}.

\begin{proposition} \lb{P.3.3}
Let $\til \ell$ be the length of the shortest interval $I$ such that
there exists a non-zero $\psi: I\to [0,1]$, vanishing at the edges
of $I$, such that inside $I$
\begin{equation} \label{3.19}
\kappa\psi''+M f(\psi) = 0.
\end{equation}
Then $\ell\ge\til \ell$.
\end{proposition}

\begin{proof}
Assume $\phi$ is a solution of \eqref{3.8} on $(x,y)\in\bbR\times
[0,l]$ with Dirichlet boundary conditions at $y=0,l$ and
$\phi(0,\cdot,\cdot)$ compactly supported (and taking values in
$[0,1]$), such that $\phi$ does not go uniformly to zero as
$t\to\infty$. Let $\til\phi$ be the solution of \eqref{3.8} with
the same boundary conditions, but with $\til\phi(0,x,y)\equiv
\sup_x\phi(0,x,y)$. By comparison theorems, $\til\phi\ge\phi$, and
so $\til\phi$ also does not go uniformly to zero as $t\to\infty$.

Moreover, obviously $\til\phi(t,x_1,y)=\til\phi(t,x_2,y)$ for any
$t,y,x_1,x_2$, and so $\til\psi(t,y)\equiv\til\phi(t,x,y)$ is
well-defined and solves
\begin{equation} \label{3.155}
\til\psi_t=\kappa\til\psi_{yy}+M f(\til\psi).
\end{equation}
Since $\til\psi$ does not go uniformly to $0$, Proposition
\ref{P.3.5} provides us $\psi$ solving \eqref{3.19}, defined on
$[0,l]$. A simple shooting argument can be used to prove that the
set of all $l$ for which solution of \eqref{3.19} does not exist is
open. Thus the set of all $l$ for which such solution exists has a
minimum $\til\ell$, and the result follows.
\end{proof}

\begin{corollary} \lb{C.3.4}
With the above notation (and $f(T)\le T$) we have
$\ell>\pi\sqrt{\ka/M}$.
\end{corollary}

\smallskip
Remark. From results in \cite{CKR} it follows that $\ell\le
c\sqrt{\ka/M}$ for some constant $c$ depending on $f$. It follows
that the critical plateau length $\ell$ is of the order of the {\it
laminar front width} $\sqrt{\ka/M}$.

\begin{proof}
By Proposition~\ref{P.3.3}, there exists a solution $\psi$ of
\eqref{3.19} on $[0,\til\ell]$ vanishing at $0,\til\ell.$
Since $f(\psi)\le \psi$, we then have
\[
\kappa\psi''+M\psi \ge 0.
\]
That is, the lowest eigenvalue of $-\triangle$ on $[0,\til\ell]$ is
at most $M\kappa$. Hence $\til\ell\ge\pi\sqrt{\ka/M}$. But if
$\til\ell=\pi\sqrt{\ka/M}$, then necessarily
$\psi(y)=c\sin(\sqrt{M/\ka}\,y)$. This contradicts \eqref{3.19}
because $f(\psi)=0$ for small $\psi$. Hence
$\til\ell>\pi\sqrt{\ka/M}$ and Proposition \ref{P.3.3} gives the
result.
\end{proof}

The above provides also a criterion for strong quenching (recall
Definition \ref{D.1.1}).

\begin{theorem} \lb{C.3.2}
If the longest plateau of $u$ is shorter than $\til\ell$ from
Proposition \ref{P.3.3}, then $u$ is strongly quenching.
\end{theorem}

\begin{proof}
Going through the proof of Proposition \ref{T.3.1a}, one observes
that $A_0$ depends on $\tau_1$, $\tau_2$, and $L$. For $L\ge 1$
one can make $\tau_1$ only depend on $\tau_2$ and from the
conditions on $A_0$ one then sees that as long as $\tau_2$ is
bounded, $A_0$ only depends on $L$. Moreover, this dependence is
linear by Lemma \ref{L.2.2} (and Lemma \ref{L.2.1}) since the
lengths of the intervals in \eqref{3.133} and \eqref{3.144} are
$O(L)$. Therefore as long as $\tau(l,L)$ is bounded in $L$, we
have $A_0\le O(L)$ and so $u$ is strongly quenching.

Assume that the longest plateau of $u$ is has length $l$ and let
$\del\equiv(\til\ell-l)/3$. Let $\til\psi$ solve \eqref{3.155} on
$[-\del,l+\del]$ with Dirichlet boundary conditions at
$y=-\del,l+\del$ and $\psi(0,y)\ge\chi_{[0,l]}(y)$. By
$l+2\del<\til\ell$, the definition of $\til\ell$, and Proposition
\ref{P.3.5}, $\til\psi$ goes uniformly to 0 and so there is
$\til\tau<\infty$ such that
$\|\til\psi(\til\tau,\cdot)\|_\infty\le\tht_0/2$. As in the proof of
Proposition \ref{P.3.3}, comparison theorems show that if
$\phi(t,x,y)$ is any solution of \eqref{3.8} on $\bbR\times[0,l]$
with Dirichlet boundary conditions at $y=0,l$ and taking values in
$[0,1]$, then $\phi(t,x,y)\le\til\psi(t,y)$. Hence $\tau(l,L)\le
\til\tau$ for all $L$, and the result follows.
\end{proof}

The following proposition relates dynamical properties of
reaction-diffusion equation with Dirichlet boundary conditions to
existence of stationary solutions. Since we were not able to find
this simple and natural result in the literature, we provide the
proof in a slightly more general setting than needed for our
application. Let $\Omega$ be a bounded domain in $\bbR^n$ with a
smooth boundary.
We also assume for the sake of simplicity that the reaction
function is smooth. In one dimension this requirement can be
removed and $f$ only continuous is sufficient. This can be done by
approximation from above with smooth $f$, comparison principles,
and a simple ODE shooting argument.

\begin{proposition} \lb{P.3.5}
Assume that there is a solution $\phi$ of
\begin{equation} \label{3.20}
\phi_t=\kappa\Delta\phi+M f(\phi)
\end{equation}
on $(x,t)\in \Omega \times \bbR^+$, with Dirichlet boundary
conditions at $\partial \Omega$ and  $\phi(\cdot,0)$ compactly
supported (and taking values in $[0,1]$), such that $\phi$ does
not go uniformly to zero as $t\to\infty$. Then there exists a
positive solution $\psi:\Omega\to[0,1]$ of
\begin{equation} \lb{3.21}
\kappa\triangle\psi+M f(\psi)=0
\end{equation}
satisfying Dirichlet boundary conditions on $\partial \Omega.$
\end{proposition}

\begin{proof}
For the sake of simplicity we let $\kappa=M=1$. Since by the maximum
principle $\phi(x,t) \leq 1$ for any $t,$ standard regularity
estimates imply that all Sobolev norms of $\phi(x,t)$ are uniformly
bounded in time: $\|\phi(x,t)\|_{H^s(\Omega)} \leq C_s.$ Define
$\phi_-(x,t) = \limsup_{t \rightarrow \infty} \phi(x,t)$ at every $x
\in \Omega.$ We claim that $\phi_-(x,t)$ is Lipshitz continuous and
is moreover a weak subsolution, that is
\[ \int_\Omega D\phi_-(x) Dv(x) \,dx \leq \int_\Omega f(\phi_-(x))v(x)\,dx \]
for any $v \in C_0^\infty(\Omega).$ To avoid certain degenerate
cases, we define here Lipshitz continuity as
$|\phi_-(x,t)-\phi_-(y,t)| \leq C|x-y|$ for any $x,y$ which belong
to some ball $B \subset \Omega,$ with the constant $C$ independent
of $x,y$ and $B.$ Indeed, let $C_1$ be a uniform upper bound on
$|\nabla \phi(x,t)|$. Assume there exist $x,y \in B \subset
\Omega$ with $|\phi_-(x) - \phi_-(y)| > 2C_1 |x-y|.$ From the
definition of $\phi_-$ it follows that there exist $t_{n}
\rightarrow \infty$ such that either $\phi(y,t_n) - \phi_-(x) >
2C_1 |x-y|$ or $\phi(x,t_n) - \phi_-(y) > 2C_1 |x-y|.$ But this
implies that for any $\epsilon >0,$ for all sufficiently large $n$
we have $|\phi(y,t_n)-\phi(x,t_n)|> 2C_1 |x-y|-\epsilon,$ which
contradicts the bound on the gradient of $\phi.$

Notice also that compactness of $\Omega$ and uniform boundedness
of $|\nabla\phi|$ show that $\phi_-$ is not identically zero and
vanishes on $\partial\Omega$.

Define $\Delta_\delta \phi_-(x) = \delta^{-2} \sum_{j=1}^n
(\phi_-(x+\delta e_j) + \phi_- (x-\delta e_j) - 2 \phi_-(x)),$
where $e_j$ are unit vectors in coordinate directions. Next, we
claim that for any $x$ such that $\dist(x, \partial \Omega) >
\delta,$ we have $\Delta_\delta \phi_-(x) \geq -f(\phi_-(x)) -
\gamma(\delta),$ where $\gamma(\delta)$ converges to zero when
$\delta$ goes to zero. Indeed, by definition of $\phi_-(x),$ we
have that for any $\epsilon >0,$ there exists a sequence $t_n
\rightarrow \infty$ such that $|\phi_-(x) - \phi(x,t_n)| <
\epsilon$ and $\phi_-(y) \geq \phi(y,t_n) -\epsilon$ for any $y.$
Moreover, we can choose $t_n$ so that $|\phi_t(x,t_n)| <
\epsilon.$ Now
\begin{align*}
 \Delta_\delta \phi_-(x) &= \delta^{-2} \sum_j (\phi_-(x+\delta e_j) + \phi_-(x-\delta e_j) -
 2\phi_-(x))
\\ & \geq
-C\epsilon \delta^{-2} + \delta^{-2} \sum_j (\phi(x+\delta
e_j,t_n) + \phi(x-\delta e_j,t_n) - 2\phi(x,t_n)).
\end{align*}
Using the mean value theorem and uniform upper bounds on
derivatives of $\phi,$ it is not hard to show that
\[ \delta^{-2}\sum_j (\phi(x+\delta e_j,t_n) + \phi(x-\delta e_j,t_n) - 2\phi(x,t_n))
\rightarrow \phi_{x_j x_j}(x,t_n) \]
uniformly in $x$ and $t_n$ as $\delta \rightarrow 0,$ with an error bounded by $C\delta.$
Therefore,
\begin{align*}
 \Delta_\delta \phi_-(x) &\geq -C\epsilon \delta^{-2} - C\delta + \Delta \phi(x,t_n)
 \\ &\geq -f(\phi(x,t_n))
-C(\epsilon \delta^{-2} +\delta) -\epsilon
\\ & \geq
-f(\phi_-(x))-C(\epsilon \delta^{-2} +\delta +\epsilon).
\end{align*}
Since $\epsilon$ is arbitrary, this leads to $\Delta_\delta \phi_-(x) \geq -f(\phi_-(x))-C\delta.$

Given $v \in C_0^\infty (\Omega),$ $v \geq 0,$ such that
$\dist({\rm supp} (v), \partial \Omega) \geq \delta,$ we have
\[ -\int_\Omega \Delta_\delta \phi_-(x) v(x)\,dx \leq \int_\Omega f(\phi_-(x))v(x)\,dx +C\delta\|v\|_{L^1(\Omega)}. \]
Carrying out discrete integration by parts on the left hand side and passing to the limit $\delta \rightarrow 0,$
we get
\[ \int_\Omega D\phi_-(x) Dv(x)\,dx \leq \int_\Omega f(\phi_-(x)) v(x)\,dx. \]
Passage to the limit is justified since we know that $\phi_-(x)$ is Lipshitz and therefore belongs to the Sobolev
space $W^{1,\infty}.$
Thus we see that $\phi_-(x)$ is a weak subsolution of
\eqref{3.21}.

Now consider initial data $\tilde\phi(x,0)$ such that $\phi_-(x)
\leq \tilde\phi(x,0) \leq 1.$ By the maximum principle, for all
$t$ we have $\tilde\phi(x,t) \geq \phi_-(x).$ Consider $\phi_+(x)
= \liminf_{t \rightarrow \infty} \tilde\phi(x,t) \geq \phi_-(x).$
By repeating the same arguments as above, we find that $\phi_+(x)$
is a weak supersolution. Then by well-known results (see e.g.
\cite{Evans}, Theorem 9.3.1), there exists a weak solution
$\psi(x)$ of \eqref{3.21}, satisfying $\phi_-(x) \leq \psi(x) \leq
\phi_+(x).$ By boundary regularity results, $\psi(x)$ is regular
on all of $\Omega.$
\end{proof}

Results in this section extend without change to the case of shear
flows in higher dimensions. The proofs are identical to those above,
this time using higher dimensional Brownian motion. Assume that
$T(t,x,y)$ is a solution of \eqref{3.1}, \eqref{3.4} on
$\bbR_0^+\times\bbR\times\bbR^n$ with the $C^1$ shear flow $u$
satisfying $u(y)=u(y+h_je_j)$ for $j=1,2,\dots,n$ and some $h_j>0$
($\{e_1,\dots,e_n\}$ being the standard basis in $\bbR^n$).

The definition of quenching flows is identical to that for $n=1$. A
plateau of $u$ is any maximal domain $\Omega$ on which $u$ is
constant. We also say that a domain $\Omega\subseteq\bbR^n$ is {\it
quenching} if for every $L<\infty$ there is $\tau(\Omega,L)<\infty$
such that any solution $\phi$ of \eqref{3.8} on
$(x,y)\in\bbR\times\Omega$ with Dirichlet boundary conditions on
$\partial\Omega$ and $\phi(0,\cdot,\cdot)$ supported in
$[-L,L]\times\Omega$ (and taking values in $[0,1]$), satisfies
$\phi(t,x,y)\le \tht_0/2$ for $t\ge\tau(\Omega,L)$. Of course, the
quenching property again depends on $M,\ka,f$. We then have

\begin{theorem} \lb{T.3.6}
With the above notation the following hold.
\begin{SL}
\item[{\rm{(i)}}] If every plateau of $u$ is
quenching, then $u$ is quenching.
\item[{\rm{(ii)}}] If $u$ has a plateau that is not
quenching, then $u$ is not quenching.
 \item[{\rm{(iii)}}] If $\partial\Omega$ and $f$ are smooth, and there is no
non-zero $\psi:\Omega\to[0,1]$ satisfying
\[
\kappa\triangle\psi+M f(\psi) = 0
\]
and vanishing on $\partial\Omega$, then the domain $\Omega$ is
quenching. Moreover, if each plateau of $u$ is contained in some
such domain, then $u$ is strongly quenching.
\end{SL}
\end{theorem}

\smallskip
{\it Remark.} Note that if $n\ge 2$, then even non-constant $u$
can have unbounded plateaux.
\smallskip

Finally we note that we only considered initial conditions
\eqref{3.4} for the sake of simplicity of presentation. It is
obvious that our results apply also in the case of smooth initial
conditions satisfying, for instance,
\[
\chi_{[-L,L]}(x)\le T(0,x,y)\le \chi_{[-L-c_L,L+c_L]}(x).
\]

If we wish to consider initial temperatures that are not maximal
(but still above the ignition temperature $\tht_0$) on an
increasing family of regions, for example,
\[
\eta\chi_{[-L,L]}(x)\le T(0,x,y)\le \eta\chi_{[-L-c_L,L+c_L]}(x).
\]
for some $\eta\in(\tht_0,1)$, then there is only one change ---
$\ell$ in Theorem \ref{T.3.1} is defined in terms of Dirichlet
solutions $\phi$ initially compactly supported and initially
bounded above by $\eta$. The above method actually applies in the
case of any family of compactly supported initial conditions
$T_L(0,x,y)$ as long as these are such that for any $L_1$ and
$\del_1>0$ there are $L_2$ and $\del_2>0$ so that
$T_{L_2}(0,x_2,y_2)\ge T_{L_1}(0,x_1,y_1)-\del_1$ whenever
$|(x_2,y_2)-(x_1,y_1)|<\del_2$ (in particular,
$T_L(0,\cdot,\cdot)$ continuous will do). This last condition is
necessary for our proof of part (i) of Theorem \ref{T.3.1} because
now we have
\[
\Phi(t,x,y) = \bbE \bigg( T_L \bigg( 0, W_{2\ka
t}^x-\frac{A}{2\ka}\int_0^{2\ka t} u(W_s^y)ds,W_{2\ka t}^y \bigg)
\bigg).
\]
Here $\ell$ is defined in terms of $\phi$ initially bounded above
by the $T_L$'s.

\section{Scaling} \lb{S.4}

In this section we study the dependence of the ``quenching
amplitude'', that is, the infimum of all $A$ such that initial
temperature distribution
\begin{equation} \lb{4.1}
T(0,x,y)=\chi_{[-L,L]}(x)
\end{equation}
leads to quenching, on the scaling of the profile of the shear
flow $u$. Hence we consider
\begin{equation} \lb{4.2}
T_t  = \ka \triangle T - Au(\al y)T_x + M f(T)
\end{equation}
with $u$ periodic and $\al>0$. The results of this section are
motivated by and agree well with numerical simulations performed
in \cite{VCKRR}. The first is

\begin{theorem} \lb{T.4.1}
Let $u\in C(\bbR)$ be a periodic function with period $h$. Then
there is $C>0$ such that for large enough $\alpha$ and $|A|\ge
C\alpha L$, the solution of \eqref{4.2} with initial condition
\eqref{4.1} satisfies $T(t,x,y)\to 0$ as $t\to \infty$, uniformly
in $\bbR^2$.
\end{theorem}

\smallskip
{\it Remark.} The necessity of this bound can be explained by the
fact that fast oscillations in the advection homogenize
propagation of the flame (w.r.t.~$y$) and so larger advection
amplitudes are needed to expose the hot region to diffusion.

\begin{proof}
Notice that we have
\begin{align}
\sup_x \Psi \bigg( \frac 1{2\ka},x,y \bigg) & = \sup_x \bbP \bigg(
x-\frac{A}{2\ka}\int_0^{1} u(\al W_s^y)ds\in [-L,L] \bigg)
\notag
\\ & = \sup_a \bbP \bigg( \int_0^{1} u(\al
W_s^y)ds\in \bigg[ a,a+\frac{4\ka L}{|A|} \bigg] \bigg) \notag
\\ & = \sup_a \bbP \bigg( \frac 1{\al^2} \int_0^{\al^2} u(
W_s^{\al y})ds\in \bigg[ a,a+\frac{4\ka L}{|A|} \bigg] \bigg)
\notag
\\ & = \sup_a \bbP \bigg( \frac 1{\al} \int_0^{\al^2} u(
W_s^{\al y})ds\in \bigg[ a,a+\frac{4\ka \al L}{|A|} \bigg] \bigg).
\notag 
\end{align}
Let us estimate the last integral.

First, we can assume $\int_0^h u(y)dy=0$, since, as before,
changing $u$ by a constant does not change the result. Second, let
$v(y)$ be such that $v'(y)=u(y)$ and $\int_0^h v(y)dy=0$, and
define $z(y)\equiv \int_0^y v(s)ds$. Hence, all three functions
are periodic with period $h$.

Now by the It\^ o formula (see, e.g., \cite[Proposition
1.1.4]{Nu}),
\begin{equation*}
z(W^y_t)-z(y)=\int_0^t v(W^y_s)dW^y_s + \frac 12 \int_0^t u(W^y_s)
ds
\end{equation*}
almost surely. Thus,
\begin{equation*}
\frac 1{\al} \int_0^{\al^2} u( W_s^{y})ds = \frac 2\al \big(
z(W^y_{\al^2})-z(y) \big) - 2\calM(y,\al,W^y)
\end{equation*}
with
\begin{equation*}
\calM(y,\al,W^y) \equiv \frac 1\al \int_0^{\al^2} v(W^y_s)dW^y_s.
\end{equation*}
Therefore with $c\equiv \|z\|_\infty$ we have
\begin{equation*}
\sup_x \Psi \bigg( \frac 1{2\ka},x,y \bigg) \le \sup_a \bbP \bigg(
\calM(\al y,\al,W^{\al y}) \in \bigg[ a,a+\frac{2\ka \al
L}{|A|}+\frac {4c}\al \bigg] \bigg).
\end{equation*}
From \eqref{3.5} and \eqref{3.6} we can see that to obtain
\eqref{3.7} for $\tau=(2\ka)^{-1}$ (and hence \eqref{1.2}), we
only need to prove
\begin{equation*}
\sup_{y,a} \bbP \bigg( \calM(y,\al,W^y) \in \bigg[ a,a+\frac
{2\ka}C+\frac {4c}\al \bigg] \bigg) \le \tht_0e^{-M/2\ka}
\end{equation*}
for some $C$ and all large enough $\al$. That is,
\begin{equation} \lb{4.4}
\sup_{y,a} \bbP ( \calM(y,\al,W^y) \in [ a,a+\eps ] ) \le
\tht_0e^{-M/2\ka}
\end{equation}
for small $\eps$ and large $\alpha$. However, for each $y,$ the
family $\al \calM(y, \al, W^y)$ is a martingale with respect to
$\alpha.$ It is not difficult to check that the central limit
theorem for martingales (see, e.g. \cite{Durrett}, Theorem 7.7.3, or
\cite{Shiryayev}) applies to $\calM(y, \al, W^y)$ giving convergence
in distribution to the normal random variable with variance
\[
\sigma^2 = \frac1h \int_0^h \bbE \left[ \int_0^1  v(W^z_s)^2 \,ds
\right] dz = \frac1h \int_0^h |v(z)|^2 \, dz>0,
\]
where $\bbE$ denotes expectation with respect to the Brownian
motion starting at $z.$ Moreover the convergence can be shown to
be uniform in $y$ since all the estimates entering the proof are
uniform in $y.$ This implies the estimate \eqref{4.4}.
\end{proof}

Next, we consider scaling in the opposite direction, that is
$\alpha\to 0$.

\begin{theorem} \lb{T.4.2}
If $u\in C^{n+1}(\bbR)$ is periodic and
$|u'(y)|+|u''(y)|+\dots+|u^{(n)}(y)|>0$ for some $n$ and all $y$,
then there is $C>0$ such that for small enough $\alpha>0$ and
$|A|\ge C\alpha^{-n} L$, the solution of \eqref{4.2} with initial
condition \eqref{4.1} satisfies $T(t,x,y)\to 0$ as $t\to \infty$,
uniformly in $\bbR^2$.
\end{theorem}

Let us give a short explanation of this result. Consider first the
situation as in \cite{VCKRR}, where $u(y)=\sin y$ was analyzed
numerically. When there is no flow, the critical quenching size,
according to results of Kanel', is of the order $\ell.$ Therefore
one expects that to quench initial data of size $L,$ the flow
should be able to thin it down to width $\ell$, given by Theorem
\ref{T.3.1}, in time $\tau \sim M^{-1}$ (before the reaction picks
up). The differential of velocities near the tip at points which
are distance $\sim \ell$ apart is $A \alpha^2 \ell^2,$ so we get
the condition for quenching $A \alpha^2 \ell^2 \tau \sim L,$ which
is consistent with our theorem (since the assumptions are
satisfied with $n=2$ when $u(y)=\sin y$). In a more general
setting, assume $u$ is smooth enough and $u'(0)=
u''(0)=\dots=u^{(n-1)}(0)=0$ (and $u'$ does not vanish to a higher
degree elsewhere). If then $A$ grows slower than $O(\alpha^{-n})$
as $\alpha\to 0$, the functions $Au(\alpha y)$ become very flat on
intervals around $0$ with increasing lengths. The heuristic
reasoning above then suggests that one should not expect quenching
for small $\alpha$'s.

To prove the theorem, we will need an auxiliary lemma. For
$b\in\calS^{n-1}$, the unit sphere in $\bbR^n$, we define
\[
P_b(y)\equiv b_ny^n+b_{n-1}y^{n-1}+\dots +b_1y.
\]

\begin{lemma} \lb{L.4.3}
Given any $t>0$ and $K<\infty$ we have
\begin{equation}\label{4.5}
\sup_{b\in\calS^{n-1},a} \bbP \bigg( \int_0^t P_b(W_s^0)ds \in
[a,a+\eps] \,\bigg|\, |W_s^0|\le K \text{ for } s\in[0,t] \bigg)
\to 0
\end{equation}
as $\eps\to 0$.
\end{lemma}

\begin{proof}
We define
\[
q(b,a,\eps)\equiv \bbP \bigg( \int_0^t P_b(W_s^0)ds \in [a,a+\eps]
\,\bigg|\, |W_s^0|\le K \text{ for } s\in[0,t] \bigg),
\]
and we let $N\equiv K^{n}+K^{n-1}+\dots+K$ so that
$|P_{b+\del}(y)-P_b(y)|\le N|\del|$ whenever $|y|\le K$. Hence we
need to show that, just as $p$ in Section \ref{S.2}, $q\to 0$ as
$\eps\to 0$, uniformly in $(b,a)\in\calS^{n-1}\times\bbR$. Notice
that we do not need to exclude the value $tP_b(0) = 0$ in the above
probability because the $P_b$'s have no plateaux.

The proof is identical to that of Lemma \ref{L.2.2}.
First, the absence of plateaux in the $P_b$'s gives $q(b,a,0)=0$.
Then with $\del\equiv tN|\del_1|+|\del_2|+|\del_3|$ we have
\begin{align*}
& |q(b+\del_1,a+\del_2,\eps+\del_3)-q(b,a,\eps)|
\\ & \le \bbP \bigg(
\int_0^t P_b(W_s^0)ds \in
[a-\del,a+\del]\cup[a+\eps-\del,a+\eps+\del] \,\bigg|\, |W_s^0|\le
K \text{ for } s\in[0,t] \bigg)
\end{align*}
which goes to zero as $\del\to 0$ because
$q(b,a,0)=q(b,a+\eps,0)=0$. Thus, $q$ is jointly continuous in
$(b,a,\eps)$. This means that $q(b,a,\eps)\to 0$ as $\eps\to 0$,
uniformly in any compact subset of $\calS^{n-1}\times\bbR$. Finally,
$q(b,a,\eps)=0$ for $|a|>tN+\eps$, finishing the proof.
\end{proof}

\begin{proof}[Proof of Theorem \ref{T.4.2}]
Since $u\in C^{n+1}(\bbR)$ and is periodic,
$|u'(y)|+|u''(y)|+\dots+|u^{(n)}(y)|>\rho$ for some $\rho>0$ and
all $y$. Let $K$ be such that
\begin{equation*}
\bbP \bigg( |W_s^0|\le K \text{ for } s\in[0,1] \bigg) \ge
1-\frac{\tht_0}2 e^{-M/2\ka}.
\end{equation*}
Let $C>0$, $|A|\ge C\al^{-n} L$, and
$c\equiv\|u^{(n+1)}\|_\infty/(n+1)!$. Then if $b_k\equiv
u^{(k)}(\al y)/k!$ for $k=1,\dots,n$, Taylor's theorem gives us
\[
u(\al (y+\del)) = u(\al y) + P_b(\al\del) + \til
c\al^{n+1}|\del|^{n+1}
\]
for some $|\til c|\le c$. Notice that $b$ need not be a unit
vector here.

With all the following probabilities conditioned by $|W_s^0|\le K$
for $s\in[0,1]$, we have
\begin{align}
\sup_{x,y} & \Psi \bigg( \frac 1{2\ka},x,y \bigg) \notag
\\ &\le \sup_{x,y}
\bbP \bigg( x-\frac{A}{2\ka}\int_0^{1} u(\al (y+W_s^0))ds\in [-L,L]
\bigg) + \frac{\tht_0}2 e^{-M/2\ka} \notag
\\ & \le \sup_{a,y} \bbP \bigg( \int_0^{1} \al^{-n}u(\al(y+
W_s^0))ds\in \bigg[ a,a+\frac {4\ka}{C} \bigg] \bigg) +
\frac{\tht_0}2 e^{-M/2\ka} \notag
\\ & \le \sup_{a,y} \bbP \bigg( \int_0^{1} \al^{-n} P_b(\al W_s^0)ds\in
\bigg[ a,a+\frac {4\ka}{C} + 2c\al K^{n+1} \bigg] \bigg) +
\frac{\tht_0}2 e^{-M/2\ka}  \notag
\\ & = \sup_{a,y} \bbP \bigg( \int_0^{1} P_d(W_s^0)ds\in
\bigg[ a,a+\frac {4\ka}{C} + N\al \bigg] \bigg) + \frac{\tht_0}2
e^{-M/2\ka} \lb{4.6}
\end{align}
with $d_k\equiv b_k\al^{k-n}$ and $N\equiv 2cK^{n+1}$. If we take
$\al<1$, then $|d|\ge |b|\ge\rho/(n+1)!$ and so there are
$e\in\calS^{n-1}$ and $r\ge\rho/(n+1)!$ such that $d=re$. The last
expression in \eqref{4.6} is then at most
\[
\sup_{e\in\calS^{n-1},a} \bbP \bigg( \int_0^{1} P_e(W_s^0)ds\in
\bigg[ a,a+ \frac{(n+1)!}\rho \bigg( \frac {4\ka}{C} + N\al \bigg)
\bigg] \bigg) + \frac{\tht_0}2 e^{-M/2\ka}.
\]
Lemma \ref{L.4.3} ensures that for some $C<\infty$ and all small
$\al$ the supremum is smaller than $\tht_0 e^{-M/2\ka}/2$, and then
\eqref{3.5} and \eqref{3.6} give \eqref{3.7} for $\tau=(2\ka)^{-1}$.
The result follows.
\end{proof}

{\bf Acknowledgement.} \rm We thank Peter Constantin, Tom Kurtz,
David Nualart, and Lenya Ryzhik for useful communications. AK has
been supported in part by NSF grants DMS-0321952 and DMS-0314129,
and Alfred P. Sloan fellowship. AZ has been supported in part by
NSF grant DMS-0314129.

\end{document}